\documentclass[10pt]{amsart}
\usepackage{amsfonts}
\usepackage{ifthen}
\usepackage{amsthm}
\usepackage{amsmath,mathrsfs}
\usepackage{graphicx}
\usepackage{amscd,amssymb,amsthm}
\usepackage{color}
\usepackage{hyperref}

\newcounter{minutes}
\divide\time by 60
\newcounter{hours}
\multiply\time by 60 \addtocounter{minutes}{-\time}

\newtheorem{theorem}{Theorem}
\newtheorem{lemma}{Lemma}
\newtheorem{corollary}{Corollary}

\newtheorem{conjecture}{Conjecture}

\title{Redheffer type bounds for Bessel and modified Bessel functions of the first kind}

\author[\'A. Baricz]{\'Arp\'ad Baricz}
\address{Institute of Applied Mathematics, \'Obuda University, 1034 Budapest, Hungary}
\address{Department of Economics, Babe\c{s}-Bolyai University, 400591 Cluj-Napoca, Romania}
\email{bariczocsi@yahoo.com}

\author[K. Mehrez]{Khaled Mehrez}
\address{D\'epartement de Math\'ematiques ISSAT Kasserine, Universit\'e de Kairouan, Tunisia}
\email{k.mehrez@yahoo.fr}



\begin{document}

\def\thefootnote{}
\footnotetext{ \texttt{File:~\jobname .tex,
          printed: \number\year-0\number\month-\number\day,
          \thehours.\ifnum\theminutes<10{0}\fi\theminutes}
} \makeatletter\def\thefootnote{\@arabic\c@footnote}\makeatother

\maketitle

\begin{abstract}
In this paper our aim is to show some new inequalities of Redheffer type for Bessel and modified Bessel
functions of the first kind. The key tools in our proofs are some classical results on the
monotonicity of quotients of differentiable functions as well as on the monotonicity of quotients of two
power series. We use also some known results on the quotients of Bessel and
modified Bessel functions of the first kind, and by using the monotonicity of the Dirichlet eta
function we prove a sharp inequality for the tangent function. At the end of the paper a conjecture
is stated, which may be of interest for further research.
\end{abstract}

\section{\bf Functional inequalities for Bessel and modified Bessel functions}
\setcounter{equation}{0}

Bessel and modified Bessel functions appear frequently in various problems of applied mathematics. Because of this their properties
worth to be studied also from the point of view of analytic inequalities. For a long list of applications concerning inequalities involving Bessel
and modified Bessel functions of the first kind we refer to the papers \cite{A}, \cite{BariczEdin}, \cite{BariczExpo} and to the references therein. Some of recent properties of Bessel and modified Bessel functions of the first kind arise naturally from the fact that they are the generalizations of the trigonometric
functions sine and cosine as well as of hyperbolic sine and cosine. Motivated by the papers \cite{A}, \cite{BW}, \cite{Khaled1}, \cite{zhuref} and \cite{ZH}, in this paper our aim is to present some new inequalities of Redheffer type for Bessel and modified Bessel
functions of the first kind. The key tools in our proofs are some classical and frequently used results on the
monotonicity of quotients of differentiable functions as well as on the monotonicity of quotients of two
power series. We use also some known results on the quotients of Bessel and
modified Bessel functions of the first kind, and by using the monotonicity of the Dirichlet eta
function we prove a sharp inequality for the tangent function. At the end of the paper a conjecture
is stated, which may be of interest for further research. The paper is organized as follows: this section contains the main results, section 2 is devoted
for proofs, while section 3 contains the conjectured inequality, which is actually proved there for a particular case.

We begin with the following Redheffer type inequality for Bessel functions of the first kind.

\begin{theorem}\label{th1}
If $\nu>-1$ and $|x|<j_{\nu,1},$ where $j_{\nu,1}$ stands for the first positive zero of the Bessel function
of the first kind $J_{\nu},$ then the following sharp exponential inequalities hold
\begin{equation}\label{ineq1}
\left(\frac{j_{\nu,1}^2-x^2}{j_{\nu,1}^2}\right)^{\alpha_{\nu}}\leq\mathcal{J}_{\nu}(x)\leq\left(\frac{j_{\nu,1}^2-x^2}{j_{\nu,1}^2}\right)^{\beta_{\nu}},
\end{equation}
with the best possible constants $\alpha_{\nu}=\frac{j_{\nu,1}^2}{4(\nu+1)}$ and $\beta_{\nu}=1.$
\end{theorem}

Taking in \eqref{ineq1} the values $\nu=\pm\frac{1}{2},$ and since $j_{\frac{1}{2},1}=\pi$ and $j_{-\frac{1}{2},1}=\frac{\pi}{2}$, we obtain the following inequalities in the next corollary. We note that the right-hand side of \eqref{ineq1} when $\nu=\frac{1}{2}$ is actually a Jordan inequality, which was studied intensively in the last ten years by many researchers. For more details see, for example, \cite{qi} and the references therein.

\begin{corollary}
The following inequalities hold
$$\left(\frac{\pi^2-x^2}{\pi^2}\right)^{\alpha_{\frac{1}{2}}}\leq\frac{\sin x}{x}\leq\left(\frac{\pi^2-x^2}{\pi^2}\right)^{\beta_{\frac{1}{2}}},$$
for all $|x|<\pi,$ where $\alpha_{\frac{1}{2}}=\frac{\pi^2}{6}$ and $\beta_{\frac{1}{2}}=1$ are the best possible constants, and
$$\left(\frac{\pi^2-4x^2}{\pi^2}\right)^{\alpha_{-\frac{1}{2}}}\leq\cos x\leq\left(\frac{\pi^2-4x^2}{\pi^2}\right)^{\beta_{-\frac{1}{2}}},$$
for all $|x|<\frac{\pi}{2},$ where $\alpha_{-\frac{1}{2}}=\frac{\pi^2}{8}$ and $\beta_{-\frac{1}{2}}=1$ are  the best possible constants.
\end{corollary}

We would like to mention here that the right-hand side of the inequality in \eqref{ineq1} is also optimal in the sense that if we
increase the power $2$ in that inequality to bigger even powers we will get weaker inequalities. In other words, it can be
seen that we have for $k\in\mathbb{N},$ $\nu>-1$ and $|x|<j_{\nu,1}$
$$\mathcal{J}_{\nu}(x)<\frac{j_{\nu,1}^2-x^2}{j_{\nu,1}^2}<\frac{j_{\nu,1}^4-x^4}{j_{\nu,1}^4}<{\dots}<\frac{j_{\nu,1}^{2k}-x^{2k}}{j_{\nu,1}^{2k}}<{\dots}.$$
Moreover, we note that if we consider the inequality
$$\mathcal{J}_{\nu}(x)<\left(\frac{j_{\nu,1}^{2k}-x^{2k}}{j_{\nu,1}^{2k}}\right)^{\xi_{\nu}},$$
then the constant $\xi_{\nu}=1$ is optimal since by using the Mittag-Leffler expansion \eqref{mittag} we get
$$\lim_{x\nearrow j_{\nu,1}}\frac{\log \mathcal{J}_{\nu}(x)}{\log\left(\frac{j_{\nu,1}^{2k}-x^{2k}}{j_{\nu,1}^{2k}}\right)}=
\lim_{x\nearrow j_{\nu,1}}\frac{\left(\log \mathcal{J}_{\nu}(x)\right)'}{\left(\log\left(\frac{j_{\nu,1}^{2k}-x^{2k}}{j_{\nu,1}^{2k}}\right)\right)'}=
\lim_{x\nearrow j_{\nu,1}}\frac{j_{\nu,1}^{2k}-x^{2k}}{2kx^{2k-1}}\sum_{n\geq1}\frac{2x}{j_{\nu,n}^2-x^2}=1.$$

A similar result to what we have in Theorem \ref{th1} is the following.

\begin{theorem}\label{th2}
If $\nu>-1$ and $|x|<j_{\nu,1},$ then the following sharp exponential inequalities hold
\begin{equation}\label{ineq2}
\left(\frac{j_{\nu,1}^2-x^2}{j_{\nu,1}^2}\right)^{\gamma_\nu}\leq\mathcal{J}_{\nu+1}(x)\leq\left
(\frac{j_{\nu,1}^2-x^2}{j_{\nu,1}^2}\right)^{\delta_\nu},
\end{equation}
with the best possible constants $\gamma_\nu=\frac{j_{\nu,1}^2}{4(\nu+2)}$ and $\delta_\nu=0.$
\end{theorem}

Similarly as above, if we take $\nu=\pm\frac{1}{2},$ we obtain the following results.

\begin{corollary}
The following inequalities hold
\begin{equation}
\left(\frac{\pi^2-4x^2}{\pi^2}\right)^{\gamma_{-\frac{1}{2}}}\leq \frac{\sin x}{x},
\end{equation}
for all $|x|<\frac{\pi}{2},$ where $\gamma_{-\frac{1}{2}}=\frac{\pi^2}{24}$
is the best possible constant, and
\begin{equation}
\left(\frac{\pi^2-x^2}{\pi^2}\right)^{\gamma_{\frac{1}{2}}}\leq 3\left(\frac{\sin x}{x^3}
-\frac{\cos x}{x^2}\right),
\end{equation}
for all $|x|<\pi,$ where $\gamma_{\frac{1}{2}}=\frac{\pi^2}{6}$ is the best possible constant.
\end{corollary}

Now, we are going to present another result which is similar to Theorems \ref{th1} and \ref{th2}.

\begin{theorem}\label{th3}
If $\nu>-1$ and $|x|<j_{\nu,1},$ then the following sharp inequalities are valid
\begin{equation}\label{ineq3}
\left(\frac{j_{\nu,1}^2}{j_{\nu,1}^2-x^2}\right)^{\epsilon_\nu}\leq
\frac{\left(\mathcal{J}_{\nu+1}(x)\right)^{\frac{\nu+2}{\nu+1}}}{\mathcal{J}_\nu(x)}\leq \left(\frac{j_{\nu,1}^2}{j_{\nu,1}^2-x^2}\right)^{\varepsilon_\nu}
\end{equation}
with the best possible constants $\epsilon_\nu=0$
 and $\varepsilon_\nu=1.$
\end{theorem}

If we take $\nu=-\frac{1}{2}$ in the above theorem, in particular we obtain the following result.

\begin{corollary}\label{cor3}
The following inequalities hold
$$\left(\frac{\pi^2}{\pi^2-4x^2}\right)^{\epsilon_{-\frac{1}{2}}}\leq \frac{\sin^3 x}{x^3\cos x}\leq\left(\frac{\pi^2}{\pi^2-4x^2}\right)^{\varepsilon_{-\frac{1}{2}}},$$
for all $|x|<\frac{\pi}{2},$ where $\epsilon_{-\frac{1}{2}}=0$ and $\varepsilon_{-\frac{1}{2}}=1$ are the best possible constants.
\end{corollary}

It is worth to mention here that the left-hand side of the inequality in Corollary \ref{cor3} is called as Lazarevi\'c inequality and its extension, that is,
the left-hand side of \eqref{ineq3} has been already considered in \cite{A}.

Now, we are going to show that the analogues of the inequalities \eqref{ineq1} and \eqref{ineqconj2} hold for the modified Bessel functions of the first kind.

\begin{theorem}\label{th5}
If $\nu>-1,$ $r>0$ and $|x|<r$, then the following inequalities
\begin{equation}\label{ineq5}
\left(\frac{r^2-x^2}{r^2}\right)^{\rho_{\nu}}\leq\mathcal{I}_{\nu}(x)\leq\left(\frac{r^2-x^2}{r^2}\right)^{\varrho_{\nu}}
\end{equation}
hold, where $\rho_{\nu}=0,$ and $\varrho_{\nu}=-\frac{r^2}{4(\nu+1)}$ are the best possible constants.
\end{theorem}

It is worth to mention that by using a similar approach as we did in the proof of \eqref{ineq5}, Zhu \cite[Theorem 1]{zhuref} proved that if $\nu>-1,$ $r>0$ and $x\in(0,r),$ then the following
Redheffer type inequality
$$\left(\frac{r^2+x^2}{r^2-x^2}\right)^{\alpha}\leq\mathcal{I}_{\nu}(x)\leq\left(\frac{r^2+x^2}{r^2-x^2}\right)^{\beta}$$
is valid with the best possible constants $\alpha=0$ and $\beta=\frac{r^2}{8(\nu+1)}.$ It can be seen that the right-hand side of \eqref{ineq5} is weaker than the right-hand side of the above inequality of Zhu. Now, choosing in \eqref{ineq5} the values $\nu=\pm\frac{1}{2}$ and $r\in\left\{\pi,\frac{\pi}{2}\right\}$ we obtain the next result.

\begin{corollary} The following inequalities hold
$$
\left(\frac{\pi^2}{\pi^2-x^2}\right)^{\varsigma_{\frac{1}{2}}}\leq\frac{\sinh x}{x}\leq\left(\frac{\pi^2}{\pi^2-x^2}\right)^{\tau_{\frac{1}{2}}},
$$
for all $|x|<\pi,$ where $\varsigma_{\frac{1}{2}}=0$ and $\tau_{\frac{1}{2}}=\frac{\pi^2}{6}$ are the best possible constants, and
$$
\left(\frac{\pi^2}{\pi^2-4x^2}\right)^{\varsigma_{-\frac{1}{2}}}\leq\cosh x\leq\left(\frac{\pi^2}{\pi^2-4x^2}\right)^{\tau_{-\frac{1}{2}}},
$$
for all $|x|<\frac{\pi}{2},$ where $\varsigma_{-\frac{1}{2}}=0$ and $\tau_{-\frac{1}{2}}=\frac{\pi^2}{8}$ are the best possible constants.
\end{corollary}

Now, we present the analogous of \eqref{ineqconj2} for modified Bessel functions.

\begin{theorem}\label{th6}
If $\nu>-1,$ $r>0$ and $|x|<r$, then the following inequalities hold
\begin{equation}\label{ineq6}
\left(\frac{r^2-x^2}{r^2}\right)^{\kappa_{\nu}}\leq \frac{\mathcal{I}_{\nu+1}(x)}{\mathcal{I}_{\nu}(x)}\leq\left(\frac{r^2-x^2}{r^2}\right)^{\lambda_{\nu}}
\end{equation}
with the best possible constants $\kappa_{\nu}=\frac{r^2}{4(\nu+1)(\nu+2)}$ and $\lambda_{\nu}=0.$ In particular, we obtain
$$
\left(\frac{r^2-x^2}{r^2}\right)^{\kappa_{-\frac{1}{2}}}\leq \frac{\tanh x}{x}\leq\left(\frac{r^2-x^2}{r^2}\right)^{\lambda_{-\frac{1}{2}}}
$$
with the best possible constants $\kappa_{-\frac{1}{2}}=\frac{\pi^2}{12}$ and $\lambda_{-\frac{1}{2}}=0.$
\end{theorem}

\section{\bf Proof of the main results}
\setcounter{equation}{0}

In the proof of the main results we will need the following two lemmas. The first lemma is about the monotonicity of two
power series, see \cite{ponn} for more details.

\begin{lemma}\label{lem1}
Let $\{a_n\}_{n\geq0}$ and $\{b_n\}_{n\geq0}$ two sequences of real numbers, and let the power series $f(x)=\sum_{n\geq0}a_{n}x^{n}$ and $g(x)=\sum_{n\geq0}b_{n}x^{n}$ be convergent for $|x|<r.$ If $b_n>0$ for $n\geq0$ and if the sequence $\left\{{a_n}/{b_n}\right\}_{n\geq0}$ is (strictly) increasing (decreasing), then the function $x\mapsto {f(x)}/{g(x)}$ is (strictly) increasing (decreasing) on $(0,r).$
\end{lemma}

The second lemma is the so-called monotone form of l'Hospital's rule, see \cite{and} for a proof.

\begin{lemma}\label{lem2}
Let $f,g:[a,b]\longrightarrow\mathbb{R}$ be two continuous functions which are differentiable on $(a,b)$. Further, let $g^{'}\neq0$ on $(a,b).$ If ${f^{\prime}}/{g^{\prime}}$ is increasing (decreasing) on $(a,b)$, then the functions $$x\mapsto\frac{f(x)-f(a)}{g(x)-g(a)}\quad \mbox{and}\quad x\mapsto \frac{f(x)-f(b)}{g(x)-g(b)}$$ are also increasing (decreasing) on $(a,b).$
\end{lemma}

Now, we are ready to present the proofs of the main results.

\begin{proof}[\bf Proof of Theorem \ref{th1}]
Since the expressions in (\ref{ineq1}) are even in $x$, in what follows we suppose without loss of
generality that $x\in(0,j_{\nu,1}).$ We consider the function $\varphi_{\nu}:(0,j_{\nu,1})\to\mathbb{R},$ defined by
$$\varphi_{\nu}(x)=\frac{\log \mathcal{J}_{\nu}(x)}{\log\left(\frac{j_{\nu,1}^2-x^2}{j_{\nu,1}^2}\right)}=\frac{f_{\nu}(x)}{g_{\nu}(x)}.$$
Since for all $x\in(-j_{\nu,1},j_{\nu,1} )$ and $\nu>-1$ we have $\mathcal{J}_{\nu}(x)>0$ (see \cite[Theorem 3]{A}), it follows that the function $\varphi_{\nu}$ is well defined. By using the differentiation formula \cite[p. 18]{wa}
\begin{equation}\label{differen}\mathcal{J}_{\nu}'(x)=-\frac{x}{2(\nu+1)}\mathcal{J}_{\nu+1}(x)\end{equation}
we get
$$\frac{f_{\nu}'(x)}{g_{\nu}'(x)}=\frac{\frac{x}{2}\frac{J_{\nu+1}(x)}{J_{\nu}(x)}}{\frac{x^2}{j_{\nu,1}^2-x^2}}=
\frac{\sum\limits_{m\geq1}\sigma_{\nu}^{(2m)}x^{2m}}{\sum\limits_{m\geq1}j_{\nu,1}^{-2m}x^{2m}},$$ where
$\sigma_{\nu}^{(2m)},$ $m\in\mathbb{N},$ is the Rayleigh function of order $2m$ defined by $\sigma_{\nu}^{(2m)}=\sum_{n\geq1}j_{\nu,n}^{-2m}$ and
we used the Kishore formula \cite{KI}
\begin{equation}\label{kishore}
\frac{x}{2}\frac{J_{\nu+1}(x)}{J_{\nu}(x)}=\sum_{m\geq1}\sigma_{\nu}^{(2m)}x^{2m},\quad |x|<j_{\nu,1}.
\end{equation}
Now, we consider the sequence $\alpha_{m,\nu}=j_{\nu,1}^{2m}\sigma_{\nu}^{(2m)},$ $m\in\mathbb{N},$ which satisfies
$$\alpha_{m+1,\nu}-\alpha_{m,\nu}=j_{\nu,1}^{2m}\sum_{n\geq1}j_{\nu,n}^{-2m}\left(\frac{j_{\nu,1}^{2}}{j_{\nu,n}^{2}}-1\right)<0$$
for all $\nu>-1$ and $m\in\mathbb{N}.$ By using Lemma \ref{lem1} we clearly have that ${f_{\nu}'}/{g_{\nu}'}$ is decreasing on $(0,j_{\nu,1}),$ and consequently the function $\varphi_{\nu}$ is also decreasing on $(0,j_{\nu,1}),$ by means of Lemma \ref{lem2}. On the other hand, by using the first Rayleigh sum $\sigma_{\nu}^{(2)}=\sum_{n\geq1}{j_{\nu,n}^{-2}}=\frac{1}{4(\nu+1)}$
and the Mittag-Leffler expansion
\begin{equation}\label{mittag}\frac{J_{\nu+1}(x)}{J_{\nu}(x)}=\sum_{n\geq1}\frac{2x}{j_{\nu,n}^2-x^2},\end{equation}
we obtain by the Bernoulli-l'Hospital's rule
$$\lim_{x\searrow 0}\varphi_{\nu}(x)=\lim_{x\rightarrow 0}\left((j_{\nu,1}^2-x^2)\frac{J_{\nu+1}(x)}{2xJ_{\nu}(x)}\right)=\frac{j_{\nu,1}^2}{4(\nu+1)},$$
and
$$\lim_{x\nearrow j_{\nu,1}}\varphi_{\nu}(x)=\lim_{x\nearrow j_{\nu,1}}\left((j_{\nu,1}^2-x^2)\frac{J_{\nu+1}(x)}{2xJ_{\nu}(x)}\right)
=\lim_{x\nearrow j_{\nu,1}}\sum_{n\geq1}\frac{j_{\nu,1}^2-x^2}{j_{\nu,n}^2-x^2}=1.$$

It is important to mention here that there is another proof of the inequalities (\ref{ineq1}). Namely, if we consider the function $\phi_{\nu}:(0,j_{\nu,1})\rightarrow \mathbb{R},$ defined by
$$\phi_{\nu}(x)=\frac{j_{\nu,1}^2}{4(\nu+1)}\log\left(\frac{j_{\nu,1}^2-x^2}{j_{\nu,1}^2}\right)-\log\mathcal{J}_{\nu}(x),$$
then taking into account the inequality \cite{BW,IF}
$$\frac{\mathcal{J}_{\nu+1}(x)}{\mathcal{J}_{\nu}(x)}<\frac{j_{\nu,1}^2}{j_{\nu,1}^2-x^2},$$
which holds for all $\nu>-1$ and $x\in(0,j_{\nu,1}),$ we obtain that
$$\phi_{\nu}'(x)=\frac{x}{2(\nu+1)}\left(\frac{\mathcal{J}_{\nu+1}(x)}{\mathcal{J}_{\nu}(x)}-\frac{j_{\nu,1}^2}{j_{\nu,1}^2-x^2}\right)<0,$$
where $\nu>-1$ and $x\in(0,j_{\nu,1}).$ Consequently, the function $\phi_{\nu}$ is decreasing on $(0,j_{\nu,1}),$ which in turn implies that $\phi_{\nu}(x)\leq \phi_{\nu}(0)=0.$ This proves the left-hand side of \eqref{ineq1}. Finally, observe that in view of the infinite product representation of $J_{\nu}(x)$ the right-hand side of \eqref{ineq1} can be written as
$$\mathcal{J}_{\nu}(x)=\prod_{n\geq 1}\left(1-\frac{x^2}{j_{\nu,n}^2}\right)<1-\frac{x^2}{j_{\nu,1}^2},$$
which is true since every factor in the above product is strictly less than 1 when $\nu>-1$ and $|x|<j_{\nu,1}.$
\end{proof}

\begin{proof}[\bf Proof of Theorem \ref{th2}]
We consider the function $\Phi_\nu:(0,j_{\nu,1})\to\mathbb{R},$ defined by
$$\Phi_\nu(x)=\frac{\log \mathcal{J}_{\nu+1}(x)}{\log\left(\frac{j_{\nu,1}^2-x^2}{j_{\nu,1}^2}\right)}=\frac{q_{\nu}(x)}
{r_{\nu}(x)}.$$
By using the differentiation formula \eqref{differen}, in view of Kishore's formula \eqref{kishore} we get
\begin{equation}
\frac{q_{\nu}'(x)}{r_{\nu}'(x)}=\frac{\frac{x}{2}\frac{J_{\nu+2}(x)}{J_{\nu+1}(x)}}
{\left(\frac{x^2}{j_{\nu,1}^2-x^2}\right)}=\frac{\sum\limits_{m\geq1}\sigma_{\nu+1}^{(2m)}x^{2m}}{\sum\limits_{m\geq1}j_{\nu,1}^{-2m}x^{2m}}.
\end{equation}
Now, we define the sequence $\beta_{m,\nu}=j_{\nu,1}^{2m}\sigma_{\nu+1}^{(2m)},$ $m\in\mathbb{N}.$ We have
$$\beta_{m+1,\nu}-\beta_{m,\nu}=j_{\nu,1}^{2m}\sum_{n\geq1} j_{\nu,n}^{-2m}\left(\frac
{j_{\nu,1}^{2}}{j_{\nu+1,n}^{2}}-1\right).$$
On the other hand, we known that (see for example \cite[p. 317]{ME}) for each $n\in\mathbb{N}$ fixed,
the function $\nu\mapsto j_{\nu,n}^2$ is increasing on $(-1,\infty).$ Consequently, the sequence $\left\{\beta_{m,\nu}\right\}_{m\geq1}$ is decreasing. Thus the function $q_{\nu}'/r_{\nu}'$ is decreasing on $(0,j_{\nu,1}),$ in view of Lemma \ref{lem1}. So,
the function $\Phi_\nu$ is also decreasing on $(0,j_{\nu,1}),$ by using Lemma \ref{lem2}. Moreover,
$$\lim_{x\searrow0}\Phi_\nu(x)=\frac{\sigma_{\nu+1}^{(2)}}{j_{\nu,1}^{-2}}=\frac{j_{\nu,1}^2}{4(\nu+2)}\ \ \ \ \mbox{and}\ \ \ \
\lim_{x\nearrow j_{\nu,1}}\Phi_\nu(x)=0.$$
\end{proof}

\begin{proof}[\bf Proof of Theorem \ref{th3}]
The left-hand side of \eqref{ineq3} is known, it can be found in \cite[Theorem 3]{A}. Thus, using this and the fact that
$$\log\left(\frac{\left(\mathcal{J}_{\nu+1}(x)\right)^{\frac{\nu+2}{\nu+1}}}
{\mathcal{J}_\nu(x)}\right)=\frac{\nu+2}{\nu+1}\log\left(\frac{\mathcal{J}_{\nu+1}(x)}
{\mathcal{J}_{\nu}(x)}\right)+\frac{1}{\nu+1}\log \mathcal{J}_{\nu}(x)$$
in view of the inequalities \eqref{ineq1} and the right-hand side of \eqref{ineqconj2} we deduce that
$$
0\leq\log\left(\frac{\left(\mathcal{J}_{\nu+1}(x)\right)^{\frac{\nu+2}{\nu+1}}}
{\mathcal{J}_\nu(x)}\right)\leq\log\left(\frac{j_{\nu,1}^2}{j_{\nu,1}^2-x^2}\right).
$$
Thus, we just need to show that the constants $\epsilon_{\nu}$ and $\varepsilon_{\nu}$ are the best possible. For this we consider the function
$\Omega_\nu:(0,j_{\nu,1})\to\mathbb{R},$ defined by
$$\Omega_\nu(x)=\frac{\log\left(\frac{\left(\mathcal{J}_{\nu+1}(x)\right)^{\frac{\nu+2}{\nu+1}}}
{\mathcal{J}_\nu(x)}\right)}{\log\left(\frac{j_{\nu,1}^2}{j_{\nu,1}^2-x^2}\right)}.$$
Using the Bernoulli-l'Hospital rule and the Mittag--Leffler expansion \eqref{mittag}, we get
\begin{equation*}
\begin{split}
\lim_{x\searrow0}\Omega_\nu(x)&=\lim_{x\searrow0}(j_{\nu,1}^2-x^2)
\left(\frac{J_{\nu+1}(x)}{2xJ_\nu(x)}-\frac{\nu+2}{\nu+1}\frac{J_{\nu+2}(x)}{2xJ_{\nu+1}(x)}\right)\\
&=\lim_{x\searrow0}(j_{\nu,1}^2-x^2)\left(\sum_{n\geq1}
\frac{1}{j_{\nu,n}^2-x^2}-\frac{\nu+2}{\nu+1}\sum_{n\geq1}\frac{1}{j_{\nu+1,n}^2-x^2}\right)=0,
\end{split}
\end{equation*}
and
\begin{equation*}
\begin{split}
\lim_{x\nearrow j_{\nu,1}}\Omega_\nu(x)&=\lim_{x\nearrow j_{\nu,1}}(j_{\nu,1}^2-x^2)\left(\frac{J_{\nu+1}(x)}{2xJ_\nu(x)}-\frac{\nu+2}{\nu+1}\frac{J_{\nu+2}(x)}{2xJ_{\nu+1}(x)}\right)\\
&=\lim_{x\nearrow j_{\nu,1}}(j_{\nu,1}^2-x^2)\left(\sum_{n\geq1}\frac{1}{j_{\nu,n}^2-x^2}-\frac{\nu+2}{\nu+1}\sum_{n\geq1}
\frac{1}{j_{\nu+1,n}^2-x^2}\right)\\
&=\lim_{x\nearrow j_{\nu,1}}\left(1+\sum_{n\geq2}\frac{j_{\nu,1}^2-x^2}{j_{\nu,n}^2-x^2}-\frac{\nu+2}{\nu+1}\sum_{n\geq1}\frac{j_{\nu,1}^2-x^2}{j_{\nu+1,n}^2-x^2}\right)=1.
\end{split}
\end{equation*}
\end{proof}

\begin{proof}[\bf Proof of Theorem \ref{th5}]
We consider the function $\Psi_{\nu}:[0,r)\to\mathbb{R},$ defined by
$$\Psi_{\nu}(x)=\log(\mathcal{I}_{\nu}(x))-\frac{r^2}{4(\nu+1)}\log\left(\frac{r^2}{r^2-x^2}\right).$$
Using the differentiation formula
$$\mathcal{I}_{\nu}'(x)=\frac{x}{2(\nu+1)}\mathcal{I}_{\nu+1}(x)$$
and the Mittag-Leffler expansion for the modified Bessel functions of first kind
$$\frac{I_{\nu+1}(x)}{I_{\nu}(x)}=\sum_{n\geq1}\frac{2x}{j_{\nu,n}^2+x^2},$$
together with the first Rayleigh sum $\sigma_{\nu}^{(2)}$ we obtain
$$\Psi_{\nu}'(x)=\frac{I_{\nu+1}(x)}{I_{\nu}(x)}-\frac{2xr^2}{4(\nu+1)(r^2-x^2)}=
\sum_{n\geq1}\frac{2x}{j_{\nu,n}^2+x^2}-\sum_{n\geq1}\frac{2xr^2}{j_{\nu,n}^2(r^2-x^2)}=
\sum_{n\geq1}\frac{-2x^3(r^2+j_{\nu,n}^2)}{j_{\nu,n}^2(j_{\nu,n}^2+x^2)(r^2-x^2)}.$$	
Therefore the function $\Psi_{\nu}$ is decreasing on $[0,r),$ and hence $\Psi_{\nu}(x)\leq \Psi_{\nu}(0)=0,$ which implies the right-hand side of \eqref{ineq5}. To prove the left-hand side of \eqref{ineq5}, from the above differentiation formula we conclude that the function $x\mapsto\mathcal{I}_{\nu}(x)$ is increasing on $[0,r)$ and hence $\mathcal{I}_{\nu}(x)\geq1.$ This yields the left-hand side of \eqref{ineq5}.

Alternatively, the inequalities in \eqref{ineq5} can be proved in the following way. We consider the function $\Gamma_{\nu}:(0,r)\to\mathbb{R},$ defined by
$$\Gamma_{\nu}(x)=\frac{u_{\nu}(x)}{v_{\nu}(x)},\ \ \  \mbox{where}\ \ \
u_{\nu}(x)=\log\mathcal{I}_{\nu}(x)\ \ \ \mbox{and}\ \ \ v_{\nu}(x)=\log\left(\frac{r^2}{r^2-x^2}\right).$$
We have that
$$\frac{u_{\nu}'(x)}{v_{\nu}'(x)}=\frac{(r^2-x^2)\mathcal{I}_{\nu}'(x)}{2x\mathcal{I}_{\nu}(x)}=\frac{P_{\nu}(x)}{2Q_{\nu}(x)},$$
where
$$P_{\nu}(x)=(r^2-x^2)\mathcal{I}_{\nu}'(x)=\frac{r^2x}{2(\nu+1)}+
\sum_{n\geq1}\left(\frac{2r^2(n+1)\Gamma(\nu+1)}{2^{2(n+1)}(n+1)!\Gamma(\nu+n+2)}-\frac{2n\Gamma(\nu+1)}{2^{2n}n!\Gamma(\nu+n+1)}\right)x^{2n+1}$$
and
$$Q_{\nu}(x)=1+\sum_{n\geq1}\frac{\Gamma(\nu+1)}{2^{2n}n!\Gamma(\nu+n+1)}x^{2n+1}.$$
Using the notation $P_{\nu}(x)=\sum_{n\geq0}a_n x^{2n+1}$ and $Q_{\nu}(x)=\sum_{n\geq0}b_n x^{2n+1}$ we obtain that the sequence $\{c_n\}_{n\geq0},$ defined by $$c_n=\frac{a_n}{b_n}=\frac{r^2}{2(\nu+n+1)}-2n,$$ is decreasing, which in turn implies that (in view of Lemma \ref{lem1}) the function $u_{\nu}'/v_{\nu}'$ is decreasing on $(0,r).$ Thus, by using Lemma \ref{lem2} we get that $\Gamma_{\nu}$ is decreasing on $(0,r),$ which implies the inequalities in \eqref{ineq5}.

It remained to show that the corresponding constants in the theorem are the best possible ones. For this we note that $\lim_{x\nearrow r}\Gamma_{\nu}(x)=\rho_{\nu}=0$ and using the Bernoulli-l'Hospital's rule we have
$$\lim_{x \searrow0}\Gamma_{\nu}(x)=\lim_{x\searrow 0}\sum_{n\geq1}\frac{2x}{j_{\nu,n}^2+x^2}\frac{r^2-x^2}{2x}=\frac{r^2}{4(\nu+1)}=\varrho_{\nu}.$$
So, indeed $\rho_{\nu}=0$ and $\varrho_{\nu}=\frac{r^2}{4(\nu+1)}$ are the best possible constants.
\end{proof}

\begin{proof}[\bf Proof of Theorem \ref{th6}]
Since the function $\nu\mapsto\mathcal{I}_{\nu}(x)$ is decreasing on $(-1,\infty)$ for all $x\in\mathbb{R},$ (see \cite[Theorem 1]{A}), in particular $\mathcal{I}_{\nu+1}(x)\leq\mathcal{I}_{\nu}(x),$ and thus the right-hand side of \eqref{ineq6} is true. Now, we prove the left-hand side of \eqref{ineq6}. By using the inequality \cite[Theorem 1]{A}
$$\mathcal{I}_{\nu+1}(x)\geq\left(\mathcal{I}_{\nu}(x)\right)^{\frac{\nu+1}{\nu+2}},$$
which holds for all $\nu>-1$ and $x\in\mathbb{R},$ in view of the right-hand side of \eqref{ineq5} we obtain that
$$\frac{\mathcal{I}_{\nu+1}(x)}{\mathcal{I}_{\nu}(x)}\geq\frac{1}{(\mathcal{I}_{\nu}(x))^{\frac{1}{\nu+2}}}
\geq\left(\frac{r^2-x^2}{r^2}\right)^{\frac{r^2}{4(\nu+1)(\nu+2)}}.$$
Alternatively, we can consider the function $\Theta_{\nu}:(0,r)\to\mathbb{R},$ defined by
$$\Theta_{\nu}(x)=\frac{\log\left(\frac{\mathcal{I}_{\nu+1}(x)}{\mathcal{I}_{\nu}(x)}\right)}{\log\left(\frac{r^2-x^2}{r^2}\right)}=\frac{s_{\nu}(x)}{t_{\nu}(x)}.$$
Since the function
$$x\mapsto \frac{s_{\nu}'(x)}{t_{\nu}'(x)}=\frac{x^2-r^2}{2x}\left(\frac{I_{\nu+2}(x)}{I_{\nu+1}(x)}-\frac{I_{\nu+1}(x)}{I_{\nu}(x)}\right)=
(x^2-r^2)\cdot \sum_{n\geq 1}\frac{j_{\nu,n}^2-j_{\nu+1,n}^2}{(j_{\nu,n}^2+x^2)(j_{\nu+1,n}^2+x^2)}$$
is decreasing as a product of two negative and increasing functions, by using the monotone form of l'Hospital's rule, that is Lemma \ref{lem2}, we obtain that $\Theta_{\nu}$ is also decreasing, and thus we have the inequalities in \eqref{ineq6}. Here we used again that for every fixed $n\in\mathbb{N}$ we have that $\nu\mapsto j_{\nu,n}^2$ is increasing on $(-1,\infty).$ By using the Bernoulli-l'Hospital rule and the Mittag-Leffler expansion we get
$$\lim_{x\searrow 0}\Theta_{\nu}(x)=
\lim_{x\searrow 0}\left(\frac{x^2-r^2}{2x}\right)\left(\sum_{n\geq1}\frac{2x}{j_{\nu+1,n}^2+x^2}-\sum_{n\geq1}\frac{2x}{j_{\nu,n}^2+x^2}\right)=\kappa_{\nu}=\frac{r^2}{4(\nu+1)(\nu+2)}.$$
On the other hand,
$$\lim_{x\nearrow r}\Theta_{\nu}(x)=\lambda_{\nu}=0.$$
These limits show that the constants $\kappa_{\nu}$ and $\lambda_{\nu}$ are the best possible.
\end{proof}

\section{\bf A conjectured inequality for Bessel functions}
\setcounter{equation}{0}

Now, consider the function $\psi_{\nu}:(0,j_{\nu,1})\to\mathbb{R},$ defined by
$$\psi_{\nu}(x)=\frac{\log\left(\frac{\mathcal{J}_{\nu+1}(x)}{\mathcal{J}_{\nu}(x)}\right)}{\log\left(\frac{j_{\nu,1}^2}{j_{\nu,1}^2-x^2}\right)}.$$
Since for all $|x|<j_{\nu,1}$ and $\nu>-1$ we have $\mathcal{J}_{\nu}(x)>0$ it follows that $\mathcal{J}_{\nu+1}(x)>0$ for each $|x|<j_{\nu+1,1},$ and since $(-j_{\nu,1},j_{\nu,1})\subset(-j_{\nu+1,1},j_{\nu+1,1})$ we conclude that the function $\psi_{\nu}$ is well defined. By using the differentiation formula
\eqref{differen} and the Kishore formula \eqref{kishore} it follows that
$$\frac{\left(\log\left(\frac{\mathcal{J}_{\nu+1}(x)}{\mathcal{J}_{\nu}(x)}\right)\right)^\prime}
{\left(\log\left(\frac{j_{\nu,1}^2}{j_{\nu,1}^2-x^2}\right)\right)^\prime}
=\frac{\frac{x\mathcal{J}_{\nu+1}(x)}{2(\nu+1)\mathcal{J}_{\nu}(x)}-\frac{x\mathcal{J}_{\nu+2}(x)}{2(\nu+2)\mathcal{J}_{\nu+1}(x)}}
{\frac{2x}{j_{\nu,1}^2-x^2}}
=\frac{\frac{x}{2}\left(\frac{J_{\nu+1}(x)}{J_{\nu}(x)}-\frac{J_{\nu+2}(x)}{J_{\nu+1}(x)}\right)}{\frac{x^2}{j_{\nu,1}^2-x^2}}
=\frac{\sum\limits_{m\geq1}\left(\sigma_{\nu}^{(2m)}-\sigma_{\nu+1}^{(2m)}\right)x^{2m}}{\sum\limits_{m\geq1}j_{\nu,1}^{-2m}x^{2m}}.$$
We believe but were unable to prove (except in the case when $\nu=-\frac{1}{2}$) that $\omega_{m,\nu}=j_{\nu,1}^{2m}\left(\sigma_{\nu}^{(2m)}-\sigma_{\nu+1}^{(2m)}\right),$ $m\in\mathbb{N},$ is increasing in $m$ for each $\nu>-1.$ If this result would be true, then in view of Lemma \ref{lem1} the above quotient (which we denote by $\omega_{\nu}(x)$) would be increasing on $(0,j_{\nu,1})$ which by using the monotone form of Bernoulli-l'Hospital's rule (that is, Lemma \ref{lem2}) would imply that the function $\psi_{\nu}$ is also increasing on $(0,j_{\nu,1}).$ However, in order to show the monotonicity of the sequence $\{\omega_{m,\nu}\}_{m\geq1}$ we were not able to show that the inequality
\begin{equation}\label{ineqconj}
j_{\nu,1}^2>\frac{\sigma_{\nu}^{(2m)}-\sigma_{\nu+1}^{(2m)}}{\sigma_{\nu}^{(2m+2)}-\sigma_{\nu+1}^{(2m+2)}}
\end{equation}
is valid for all $\nu>-1$ and $m\in\mathbb{N}.$ Observe that by using the Bernoulli-l'Hospital rule we have that
$$\lim_{x\searrow 0}\psi_{\nu}(x)=\lim_{x\searrow 0}\frac{j_{\nu,1}^2-x^2}{4\mathcal{J}_{\nu}(x)\mathcal{J}_{\nu+1}(x)}\left(\frac{\mathcal{J}_{\nu}^2(x)}{\nu+1}-
\frac{\mathcal{J}_{\nu}(x)\mathcal{J}_{\nu+2}(x)}{\nu+2}\right)=\frac{j_{\nu,1}^2}{4(\nu+1)(\nu+2)}$$
and
$$\lim_{x\nearrow j_{\nu,1}}\psi_{\nu}(x)=\lim_{x\nearrow j_{\nu,1}}\frac{-\log\mathcal{J}_{\nu}(x)}{\log\left(\frac{j_{\nu,1}^2}{j_{\nu,1}^2-x^2}\right)}=\lim_{x\nearrow j_{\nu,1}}\frac{(j_{\nu,1}^2-x^2)J_{\nu+1}(x)}{2xJ_{\nu}(x)}=\lim_{x\nearrow j_{\nu,1}}\sum_{n\geq1}\frac{j_{\nu,1}^2-x^2}{j_{\nu,n}^2-x^2}=1.$$
These limits can be obtained also by noticing the fact that the above quotient $\omega_{\nu}(x)$ can be rewritten as
$$\sum_{n\geq1}\frac{(j_{\nu+1,n}^2-j_{\nu,n}^2)(j_{\nu,1}^2-x^2)}{(j_{\nu,n}^2-x^2)(j_{\nu+1,n}^2-x^2)}\quad \mbox{or}\quad
\sum_{n\geq1}\frac{j_{\nu,1}^2-x^2}{j_{\nu,n}^2-x^2}\sum_{n\geq1}\frac{4j_{\nu+1,n}^2}{(x^2-j_{\nu+1,n}^2)^2}.$$
The first representation follows simply from Mittag-Leffler expansion \eqref{mittag}, while the second representation from the formula of Skoovgard \cite{skov} on the Tur\'anian of Bessel functions of the first kind
$$1-\frac{J_{\nu}(x)J_{\nu+2}(x)}{J_{\nu+1}^2(x)}=\sum_{n\geq1}\frac{4j_{\nu+1,n}^2}{(x^2-j_{\nu+1,n}^2)^2}.$$ Thus, if the inequality \eqref{ineqconj} would be true we would get the following result of which left-hand side we state as a conjecture:

\begin{conjecture}
If $\nu>-1$ and $|x|<j_{\nu,1},$ then the following inequalities are valid
\begin{equation}\label{ineqconj2}
\left(\frac{j_{\nu,1}^2}{j_{\nu,1}^2-x^2}\right)^{\vartheta_{\nu}}\leq\frac{\mathcal{J}_{\nu+1}(x)}{\mathcal{J}_{\nu}(x)}
\leq\left(\frac{j_{\nu,1}^2}{j_{\nu,1}^2-x^2}\right)^{\eta_{\nu}},
\end{equation}
where $\vartheta_{\nu}=\frac{j_{\nu,1}^2}{4(\nu+1)(\nu+2)}$ and $\eta_{\nu}=1$ are the best possible constants.
\end{conjecture}

The above results would reduce to the following:

\begin{theorem}
For $|x|<\frac{\pi}{2}$ the following inequality
\begin{equation}\label{ineqconj3}
\left(\frac{\pi^2}{\pi^2-4x^2}\right)^{\vartheta_{-\frac{1}{2}}}\leq\frac{\tan x}{x}\leq\left(\frac{\pi^2}{\pi^2-4x^2}\right)^{\eta_{-\frac{1}{2}}},
\end{equation}
holds, where $\vartheta_{-\frac{1}{2}}=\frac{\pi^2}{12}$ and $\eta_{-\frac{1}{2}}=1$ are the best possible constants.
\end{theorem}

It is worth mentioning that the right-hand sides of \eqref{ineqconj2} and \eqref{ineqconj3} are not new. The right-hand side of \eqref{ineqconj2} was deduced by Ifantis and Siafarikas \cite{IF}, and by using a completely different method by Baricz and Wu \cite{BW}. Now, we are going to show that the sequence $\omega_{m,-\frac{1}{2}}$ is indeed increasing, which - following the above argument - implies that the left-hand side (and also the right-hand side) of \eqref{ineqconj3} is indeed true. According to Kishore \cite{KI} if $B_m$ denotes the $m$th Bernoulli number in the
even suffix notation and $G_m=2(1-2^m)B_m,$ then
$$\sigma_{-\frac{1}{2}}^{(2m)}=(-1)^m\frac{2^{2m-2}}{(2m)!}G_{2m}\ \ \mbox{and}\ \ \sigma_{\frac{1}{2}}^{(2m)}=(-1)^{m-1}\frac{2^{2m-1}}{(2m)!}B_{2m}.$$
In view of these expressions and
$$B_{2m}=(-1)^{m+1}\frac{2(2m)!}{(2\pi)^{2m}}\zeta(2m)$$
we get
\begin{align*}
\beta_{m+1,-\frac{1}{2}}-\beta_{m,-\frac{1}{2}}&=\left(\frac{\pi}{2}\right)^{2m+2}\frac{(-1)^{m+1}2^{2m}}{(2m+2)!}G_{2m+2}-
\left(\frac{\pi}{2}\right)^{2m}\frac{(-1)^{m}2^{2m-2}}{(2m)!}G_{2m}\\
&+\left(\frac{\pi}{2}\right)^{2m}\frac{(-1)^{m-1}2^{2m-1}}{(2m)!}B_{2m}-\left(\frac{\pi}{2}\right)^{2m+2}\frac{(-1)^{m}2^{2m+1}}{(2m+2)!}B_{2m+2}\\
&=\frac{(-1)^{m+1}\pi^{2m}}{2(2m)!}\left((1-2^{2m-1})B_{2m}+\frac{\pi^2(1-2^{2m+1})}{(2m+1)(2m+2)}B_{2m+2}\right)\\
&=\frac{1}{2}\left(1-2^{1-(2m+2)}\right)\zeta(2m+2)-\frac{1}{2}\left(1-2^{1-2m}\right)\zeta(2m)\\
&=\frac{1}{2}\left(\eta(2m+2)-\eta(2m)\right)>0,
\end{align*}
where $m\in\mathbb{N},$ $\zeta$ stands for the Riemann zeta function, $\eta$ stands for the Dirichlet eta function (or alternating zeta function) and can be
written as $\eta(s)=(1-2^{1-s})\zeta(s).$ Here we used the fact that the Dirichlet eta function is increasing on $(0,\infty),$ according to van de Lune \cite{lune}.

We note that following the proof of Theorem \ref{th3} it would be possible to give a short and elegant proof of the left-hand side of \eqref{ineq3} by using the left-hand side of the proved result \eqref{ineq1} and the left-hand side of the conjectured result \eqref{ineqconj2}.

Finally, it is worth to mention that the conjectured left-hand side of the inequality \eqref{ineqconj2} improves a known result from the literature. Namely, Baricz and Wu \cite{BW} proved recently the following inequalities
\begin{equation}\label{bw1}
\frac{\mathcal{J}_{\nu+1}(x)}{\mathcal{J}_{\nu}(x)}\geq\left(\frac{j_{\nu,1}^2+x^2}{j_{\nu,1}^2-x^2}\right)^{\frac{j_{\nu,1}^2}{8(\nu+1)(\nu+2)}},
\end{equation}
where $|x|<j_{\nu,1},$ $\nu\geq-\frac{7}{8},$ and $\frac{j_{\nu,1}^2}{8(\nu+1)(\nu+2)}$ is the best possible constants. The left-hand side of the inequality \eqref{ineqconj2} is better than (\ref{bw1}). This is justified by the following inequality
$$\vartheta_{\nu}\left(\log\frac{j_{\nu,1}^4}{(j_{\nu,1}^2-x^2)^2}-\log\frac{j_{\nu,1}^2+x^2}{j_{\nu,1}^2-x^2}\right)>0,$$
which holds for all $\nu>-1$ and $|x|<j_{\nu,1}.$

\end{document}